\documentclass[12pt]{amsart}
\usepackage[]{amsmath, amsthm, amsfonts, verbatim, amssymb}
\newcommand\C{{\mathbb C}}

\newcommand\dee{\partial}
\renewcommand\O{\Omega}
\newcommand\Obar{\overline{\Omega}}
\newcommand\Ot{\widetilde\Omega}

\numberwithin{equation}{section}

\begin{document}

\title[Green's and Ahlfors]
{The Green's function and the Ahlfors map}
\author[Bell]
{Steven R. Bell}

\address[]{Mathematics Department, Purdue University, West Lafayette,
IN  47907}
\email{bell@math.purdue.edu}

\subjclass{30C35}
\keywords{Szeg\H o kernel}

\begin{abstract}
The classical Green's function associated to a simply connected domain in the
complex plane is easily expressed in terms of a Riemann mapping function.  The
purpose of this paper is to express the Green's function of a finitely connected
domain in the plane in terms of a single Ahlfors mapping of the domain, which is
a proper holomorphic mapping of the domain onto the unit disc that is the analogue
of the Riemann map in the multiply connected setting.
\end{abstract}

\maketitle

\theoremstyle{plain}

\newtheorem {thm}{Theorem}[section]
\newtheorem {lem}[thm]{Lemma}

\hyphenation{bi-hol-o-mor-phic}
\hyphenation{hol-o-mor-phic}

\section{Introduction}
\label{sec1}
If one knows a Riemann map $f$ associated to a simply connected domain $\O\ne\C$
in the complex plane, then the classical Green's function associated to $\O$ is
given by
\begin{equation}
G(z,w)=-\ln\left|\frac{f(z)-f(w)}{1-\overline{f(w)}\, f(z)}\right|.
\label{eqnA}
\end{equation}
The function on the right hand side of equation~(\ref{eqnA}) can be realized as
the real part of an antiderivative of a rational function composed with
$(f(z),f(w),\overline{f(w)}\,)$.  To see this, let
$$\phi(z):= \ln \left|\frac{z-w}{1-\bar w\, z}\right|^2,$$
and notice that
$$\frac{\dee\phi}{\dee z}=
\frac{1-|w|^2}{(z-w)(1-\bar w\,z)}.$$
Define a rational function $R$ of
three complex variables via
$$R(z,w,v)=-\frac{1-wv}{2(z-w)(1-vz)},$$
and let $\alpha(z,w,\bar w)$ denote a local antiderivative of $R(z,w,\bar w)$ in the $z$ variable
defined via
$$\alpha(z,w,\bar w)=\int_{\gamma_z} R(\zeta,w,\bar w)\ d\zeta,$$
where, for $z$ and $w$ in the unit disc with $z\ne w$, $\gamma_z$ is any curve that
starts at $1$, enters the unit disc, avoids $w$, and ends at $z$.  Since $\phi$ is
real valued and vanishes on the boundary, it follows that
$$\phi(z)=\phi(z)-\phi(1)=\text{Re }2\int_{\gamma_z}\frac{\dee\phi}{\dee z}\ dz.$$
Hence, the Green's function
of $\O$ is given by
\begin{equation}
\label{eqnB}
G(z,w)=\text{Re }\alpha(f(z),f(w),\overline{f(w)}\,).
\end{equation}
Note that, although $\alpha$ is a multivalued function, its real part is single valued
on the unit disc.  (We shall give a more careful explanation of these ideas in \S2
when we generalize them to the multiply connected setting.)

On a bounded multiply connected domain, an Ahlfors mapping $f$ associated to a point $a$
in the domain is the solution to the same extremal problem that determines the Riemann
map on a simply connected domain, i.e., among all holomorphic functions mapping the
domain into the unit disc, the Ahlfors map is the unique function such that $f'(a)$ is
real and as large as possible.  We shall give more details about the Ahlfors map in \S2.
(See \cite{B1} for a construction of the Ahlfors map in the spirit of the present work.)

In this paper, we address the philosophical question, ``If one knows an Ahlfors
mapping of a multiply connected domain onto the unit disc,
then does one know the Green's function of the domain?''  In fact, we show
that if one knows an Ahlfors map $f$ and {\it finitely many complex numbers},
then one knows the Green's function of the domain.  We shall show
that the Green's function can be expressed in terms of a finite sum of the
real parts of finitely many functions which are antiderivatives of {\it algebraic\/}
functions composed with the Ahlfors map.  The finitely many basic functions
involved will be defined similarly to the function $\text{Re }\alpha(f(z),f(w),\overline{f(w)}\,)$
appearing in equation~(\ref{eqnB}), where the rational function $R$ of three variables
will be replaced by an algebraic function of three variables (see Theorem~\ref{thm1}
below).

When viewed in the correct light, our formula is seen to replace the rational function
inside equation~(\ref{eqnA}) with an algebraic function (which, unfortunately, may
depend on the domain under study).

To describe our main result, we must define three types of abelian functions.
Given an algebraic function $A(z,w,v)$ of three complex variables, we
shall say that $\alpha$ is a {\it Green antiderivative of type I\/} if
$$\alpha(z,w,\bar w)=\int_{\gamma_z} A(\zeta,w,\bar w)\ d\zeta,$$
where, for $z$ and $w$ in the unit disc, $\gamma_z$ is any curve that starts at $1$,
enters the unit disc, avoids finitely many points (which will be
specified), and ends at $z$.  In what follows, it will be known that, for fixed $w$ in the
unit disc, the algebraic function $A(\zeta,w,\bar w)$ of $\zeta$ is well defined and
analytic on a neighborhood of the unit circle and analytically continues nicely along
curves like $\gamma_z$ so that the integral makes sense.

We shall say that $\alpha$ is a {\it Green antiderivative of
type II\/} if
$$\alpha(w,\bar w)=\int_{\gamma} A(\zeta,w,\bar w)\ d\zeta,$$
where $\gamma$ is a fixed curve that starts at $1$,
stays in the closed unit disc, avoids finitely many points (which will be
specified), and terminates at another boundary point.  It will be known
that, for fixed $w$ in the unit disc, the algebraic function $A(\zeta,w,\bar w)$
of $\zeta$ is well defined and analytic on a neighborhood of the unit circle and
analytically continues nicely along curves like $\gamma_z$.

We shall say that $\alpha$ is a {\it Green antiderivative of type III\/} if it is
the antiderivative of an algebraic function which we may express via
$$\alpha(z)=\int_{\gamma_z} A(\zeta)\ d\zeta,$$
where $A(\zeta)$ is algebraic and, for $z$ in the unit disc, $\gamma_z$ is any curve
that starts at $1$, enters the unit disc, avoids finitely many points (which will be
specified), and ends at $z$.  It will be known that $\alpha$ is well defined and
analytic on a neighborhood of the unit circle and analytically continues along
curves like $\gamma_z$.

We can now state our main theorems.

\begin{thm}
\label{thm1}
Suppose that $\O$ is an $n$-connected domain in the plane ($n>1$) such that no
boundary component is a point and suppose that $f$ is an Ahlfors mapping
of $\O$ onto the unit disc.  There is a Green antiderivative $\alpha(z,w,\bar w)$
of type I, and Green antiderivatives $\alpha_j(w,\bar w)$ of type II, and Green
antiderivatives $\beta_j(z)$ of type III such that the Green's function associated to
$\O$ is given by
$$G(z,w)=\text{Re }\alpha(f(z),f(w),\overline{f(w)}\,) +
\sum_{j=1}^{n-1} \left(\text{Re }\alpha_j(f(w),\overline{f(w)}\,)\right)
\left(\text{Re }\beta_j(f(z))\right).$$
\end{thm}

This theorem reveals that if we know an Ahlfors mapping to the disc
and the finitely many coefficients that define the algebraic functions in the
background, then we can determine the Green's function of the domain.  It also
gives us a recipe for extending the Green's function past the boundary whenever
an Ahlfors map extends.  Furthermore, it shows that the Green's function can
be ``zipped'' down to a very small data set, namely, finitely many coefficients
plus the boundary values of a single Ahlfors map.  Since Ahlfors maps are
easily computed numerically (see \cite{B0} or Chapter~26 of \cite{B1}), this
theorem may lead to numerical recipes for the Green's function.

Following in the tradition of Bergman, Garabedian, Grunsky, and Schiffer
(see \cite{GS}), we have shown in \cite{B3,B3a} that many of the classical
kernel functions of complex analysis can be expressed in terms of {\it two\/}
Ahlfors mappings.  However, until now, the Green's function has eluded our
assaults upon it.  We got the impression from reading \cite{GS} that people
were looking for the results of the present paper, but they were missing the
algebraic functions that arise here from the connection between the Green's
function, the Szeg\H o kernel, and the field of meromorphic functions on the
double.

Before we can state our next result, we must make some definitions.
As is customary, let $\frac{\dee}{\dee z}$ denote the differential operator
$\frac{1}{2}\left(\frac{\dee}{\dee x}-i\frac{\dee}{\dee y}\right)$ and let
$\frac{\dee}{\dee\bar z}$ denote 
$\frac{1}{2}\left(\frac{\dee}{\dee x}+i\frac{\dee}{\dee y}\right)$.
Assume for the moment that $\O$ is a bounded domain bounded by $n$ non-intersecting
$C^\infty$-smooth Jordan curves.
Let $\gamma_j$, $j=1,\dots n$ denote the boundary curves of $\O$, where $\gamma_n$
is the outer boundary.  The functions $\omega_j$ are the standard harmonic
measure functions which are harmonic functions on $\O$ with boundary values given
by one on $\gamma_j$ and zero on $\gamma_k$ with $k\ne j$.  Let $F_j'$ denote
the holomorphic function given by $2(\dee/\dee z)\omega_j$.  It is well known
that the vector space $\mathcal F'$ of complex linear combinations of all the
$F_j'$ is $n-1$ dimensional and that $\{F_j'\}_{j=1}^{n-1}$ forms a basis for
$\mathcal F'$.  It is easy to show that a basis
$\{u_j\}_{j=1}^{n-1}$ for $\mathcal F'$ can be found such that
$$\delta_{kj}=\int_{\gamma_k} u_j(w)\ dw,$$
for $j=1,\dots,n-1$ and $k=1,\dots,n-1$
where $\delta_{kj}$ denotes the Kronecker delta.  Indeed, if we set
$u_j=\sum_{m=1}^{n-1}\sigma_{jm}F_m'$, then we need
$$\sum_{m=1}^{n-1}\sigma_{jm}\int_{\gamma_k}F_m'\ dz=\delta_{kj},$$
and it is well known that the matrix of periods $\int_{\gamma_k}F_m'\ dz$
is non-singular.  Since the outward normal derivative of $\omega_k$ is equal
to $-iF_k'(z)T(z)$ at a boundary point $z$, where $T(z)$ is a complex number
representing the unit tangent vector at $z\in b\O$ pointing in the direction
of the standard orientation, it follows that $i[\sigma_{jm}]$ is a matrix of
real coefficients, a fact that we shall need later on.  (See \cite{Kh} for a
nice presentation of these ideas and \cite[p.~11]{B2a} where the functions
$u_j$ arise in the context of the Green's function.)

Let $S(z,w)$ and $L(z,w)$ denote the Szeg\H o and Garabedian kernels associated
to $\O$, respectively.  Let $\mu_j(z)=
i\sum_{k=1}^{n-1}\sigma_{jk}\omega_k$ so that $2(\dee/\dee z)\mu_j=iu_j$.
Notice that $\mu_j$ is a real valued function on $\Obar$ because the coefficients
$i\sigma_{jk}$ are real.
Finally, let
$$\lambda_j(w)=\int_{\zeta\in\gamma_j}\frac{|S(w,\zeta)|^2}{S(w,w)}\ ds,$$
where $ds$ denotes the arc length measure.
It is proved in \cite{B2} that $\lambda_j$ is in $C^\infty(\Obar)$ and
has the same boundary values as $\omega_j$.

When $\O$ is a finitely connected domain such that no boundary component is
a point, we may define the functions $\omega_j$, $F_j'$, and the Szeg\H o
and Garabedian kernels via a biholomorphic mapping $\Phi$ which maps to a
bounded domain $\Ot$ with real analytic boundary curves as explained in
\cite{B2a}.  The transformation rule for the Szeg\H o kernel under
biholomorphic mappings yields that the functions $\lambda_j$ associated
to a smooth domain should transform via the same formula as $\omega_j$,
i.e., $\lambda_j(z)= \tilde{\lambda}_j(\Phi(z))$, where it is understood
that $\Phi$ takes the boundary curve $\gamma_j$ to the corresponding
boundary curve $\tilde{\gamma}_j$ of $\Ot$.  When $\O$ does not have 
smooth boundary, we take this to be the definition of $\lambda_j$.  With
these definitions, it was shown in \cite{B2} and \cite{B2a} that the
formulas we shall use in this work from \cite{B2} and \cite{B2a} are valid.

Our next theorem shows that the Green's function is a rather simple combination
of functions of {\it one complex variable}.  The function $\alpha$ appearing
in the theorem is the same as that appearing in Theorem~\ref{thm1}.

\begin{thm}
\label{thm2}
Suppose that $\O$ is an $n$-connected domain in the plane ($n>1$) such that no
boundary component is a point and suppose that $f$ is an Ahlfors mapping
of $\O$ onto the unit disc.  There is a Green antiderivative $\alpha(z,w,\bar w)$
of type I such that the Green's function associated to $\O$ is given by
$$G(z,w)=\text{Re }\alpha(f(z),f(w),\overline{f(w)}\,) +
\pi\sum_{j=1}^{n-1}(\omega_j(w)-\lambda_j(w))\mu_j(z).$$
\end{thm}

We shall also see along the way that the functions $\omega_j(w)$, $\lambda_j(w)$, and
$\mu_j(z)$ appearing in the formula in Theorem~\ref{thm2} are special functions
that are related to Green antiderivatives of various sorts.  We collect these
results in the following theorem.

\begin{thm}
\label{thm3}
Suppose that $\O$ is an $n$-connected domain in the plane ($n>1$) such that no
boundary component is a point and suppose that $f$ is an Ahlfors mapping
of $\O$ onto the unit disc.  The functions $\omega_j(w)$ and $\mu_j(w)$ are equal to
$\text{Re }\beta_j(f(w))$ and $\text{Re }\tau_j(f(w))$, respectively, 
where $\beta_j$ and $\tau_j$ are Green antiderivatives of
type~III.  The functions $\lambda_j(w)$ are equal to
$\text{Re }\psi_j(f(w),\overline{f(w)})$ where $\psi_j$ is a Green
antiderivative of type~II.
\end{thm}

Before we turn to the proofs of the theorems, we take some time to motivate the
proof and give a more precise description of the special functions
of types I, II, and III in the statement of the theorems.

\section{Motivation of the Proof}
\label{sec2}
As is typical in many arguments in conformal mapping theory, the subject of
this paper allows us to assume our domain $\Omega$ is a bounded domain bounded
by $n$ non-intersecting {\it real analytic\/} Jordan curves.  Indeed, an Ahlfors
mapping composed with a biholomorphic mapping is a constant times an Ahlfors
mapping.  Furthermore, the simple transformation formulas under biholomorphic mappings
for the Green's functions and all the other functions appearing in Theorems~\ref{thm1},
\ref{thm2}, and \ref{thm3} allow us to reduce to the case of real analytic boundary
by the standard change of variables.

Let $G(z,w)$ denote the classical Green's function
associated to $\O$ (with singularity $-\ln|z-w|$).  Given a point $a\in\O$, let
$f_a$ denote the Ahlfors map associated to the pair $(\O,a)$.  It is an
$n$-to-one mapping (counting
multiplicities), it extends holomorphically past the boundary of $\O$, and it
maps each boundary curve one-to-one onto the unit circle.  Furthermore,
$f_a(a)=0$, and $f_a$ is the unique function mapping $\O$ into the unit disc
maximizing the quantity $|f_a'(a)|$ with $f_a'(a)>0$.  An Ahlfors map should
be thought of as the replacement for the Riemann map in the non-simply
connected setting.  It shares all the properties of a Riemann map but one.
Instead of being one-to-one, it is $m$-to-one for the smallest possible $m$.
An Ahlfors mapping is an example of a proper holomorphic mapping to the unit
disc, meaning that, given a compact subset $K$ of the
unit disc, $f_a^{-1}(K)$ must be a compact subset of the domain.  (When the
boundary of $\O$ consists of finitely many non-intersecting Jordan curves,
then $f_a$ extends continuously to the boundary, and the properness condition is
equivalent to the condition that $f_a$ maps the boundary of $\O$ into
the boundary of the unit disc.  See \cite{R} for more of the basic facts about
proper holomorphic mappings.)

We shall
use the shorthand notation $G_z(z,w)$ and $G_{\bar z}(z,w)$ to denote
$\frac{\dee G}{\dee z}(z,w)$ and $\frac{\dee G}{\dee\bar z}(z,w)$, respectively.

If $h$ is a real valued harmonic function that extends smoothly up to the boundary of
$\O$ and is a constant $c_j$ on the boundary curve $\gamma_j$, then we may differentiate
the identity $c_j\equiv h(z(t))$ with respect to $t$ where  $z(t)$ parameterizes
$\gamma_j$ in the standard sense to see that
$$h_z(z(t))z'(t)+h_{\bar z}(z)\overline{z'(t)}=0.$$
Divide this equation by $|z'(t)|$ to obtain
$$h_z(z)T(z)+h_{\bar z}(z)\overline{T(z)}=0$$
for $z\in b\O$ (where we have taken $z=z(t)$ and $T(z)=z'(t)/|z'(t)|$ to be consistent
with our previous use of the complex unit tangent function $T$).
Apply this idea to the Green's function to see that
$$G_z(z,w)T(z)=-G_{\bar z}(z,w)\overline{T(z)}\quad\text{for }z\in b\O, w\in\O.$$
Note that $G_{\bar z}=\overline{G_z}$ so that
\begin{equation}
\label{eqnC}
G_z(z,w)T(z)=-\overline{G_{\bar z}(z,w)}\ \overline{T(z)}\quad\text{for }z\in b\O, w\in\O.
\end{equation}

Suppose that $f$ is an Ahlfors map associated to a point in $\O$.
The boundary identity
$0\equiv\ln|f|^2$ can be used just as we did for the Green's function above
to see that $f$ satisfies
\begin{equation}
\label{eqnD}
\frac{f'(z)}{f(z)}T(z)=
-\left(\overline{f'(z)}/\,\overline{f(z)}\right)\overline{T(z)}
\quad\text{for }z\in b\O.
\end{equation}
If we divide equation~(\ref{eqnC}) by equation~(\ref{eqnD}) and note that
$f(z)=1/\,\overline{f(z)}$ on the boundary (because $|f(z)|=1$ there), we see
that, for fixed $w\in\O$, the meromorphic function of $z$
$$\frac{G_z(z,w)}{f'(z)}$$
is equal to the complex conjugate of the meromorphic function
$$\frac{G_z(z,w)f(z)^2}{f'(z)}$$
on the boundary, and therefore it
extends to the double of $\O$ in the $z$ variable as a meromorphic
function.

We proved in \cite{B3} that, given a proper holomorphic
map $f$, there is another proper holomorphic mapping $F$ (which can be
taken to be an Ahlfors map) such that
$f$ and $F$ extend to the double and generate all the meromorphic
functions on the double, i.e., they form a primitive pair for the double.
Any two such functions are algebraicly dependent (see Farkas and Kra
\cite{FK}), so there is an
irreducible polynomial $P$ such that $P(f,F)\equiv0$.  This shows
that $F$ is an algebraic function of $f$.
Since $\frac{G_z(z,w)}{f'(z)}$ extends to the double, it
is a rational function of $f$ and $F$, and hence, it is an
algebraic function $A$ of $f$ alone.  This shows that
$$G_z(z,w)=f'(z)A(f(z))$$
for this fixed $w$.  Suppose $z_0$ is a point in the boundary of $\O$
such that $f(z_0)=1$.  Note that $f'$ is non-vanishing on the boundary of
$\O$.  Let $Z=\{z\in\Obar:f'(z)=0\}$ denote the branch locus of $f$
(which is a finite set) and
let $W=\{f(z):z\in Z\}$ denote the image of the branch locus.
Although $A$ is multivalued, there is a germ (or function element)
of $A(\zeta)$ at $\zeta=1$ which is holomorphic near $1$ and such that
$G_z=f'(A\circ f)$ continues to a neighborhood of $\Obar$ as a single
valued holomorphic function with a simple pole in $z$ at the point $w$
and removable singularities at each of the finitely many points in $Z$
different from $w$.  It will profit us here to consider
the various branches of $A$ that arise in this manner a little more
carefully.  Let $L$ denote the union of the set of closed line
segments in the unit disc joining each point in $W$ to the origin.
Note that $f$ is an unbranched covering of $\O\setminus f^{-1}(L)$ onto
$D_1(0)\setminus L$.  Since $f$ is $n$-to-one (counting multiplicities)
from $\O$ to $D_1(0)$, and since $f$ maps each boundary curve of $\O$
one-to-one onto the unit circle, we deduce that the set $f^{-1}(L)$ divides
$\O$ into $n$~two-connected components, one for each boundary curve.
Let $\mathcal O_k$ denote the component
with the boundary curve $\gamma_k$ as part of its boundary.  Let
$F_k=f^{-1}$ denote the holomorphic inverse of $f$ defined on
$D_1(0)\setminus L$ which maps this set biholomorphically to $\mathcal O_k$.
There are $n$ (single valued) holomorphic functions $A_1,A_2,\dots,A_n$ on
$D_1(0)\setminus L$ minus $f(w)$ which are branches of $A$ such that
$G_z(z,w)=f'(z)A_k(f(z))$ on $\mathcal O_k$.  Notice that each $A_k$
extends continuously up to $L\setminus W$ minus $f(w)$ from one-sided
neighborhoods, but that $A_k$ might have algebraic singularities
which might tend to infinity at points in $W\cup\{f(w)\}$.

If $\varphi$ is real-valued in a neighborhood of
a curve $\nu$ that starts at $a$ and ends at $b$, then
$$2\int_\nu \varphi_z\ dz=
\left(\int_\nu \varphi_x\,dx+\varphi_y\,dy\right)+i
\left(\int_\nu \varphi_x\,dy-\varphi_y\,dx\right),$$
and so
$$\varphi(b)-\varphi(a)=
2\text{Re }\left(\int_\nu \varphi_z\ dz\right).$$
Assume for the moment that $z$ is a point in $\O$ that is close to
the boundary point $z_0$ that we chose above so that $f(z_0)=1$.
Since the Green's function is real
valued and vanishes on the boundary, we may use this idea to
antidifferentiate via
$$G(z,w)=2\text{Re }\left(\int_\Gamma G_\zeta(\zeta,w)\ d\zeta\right)
=\text{Re }\left(\int_\Gamma 2f'(\zeta)A(f(\zeta))\ d\zeta\right),$$
where $\Gamma$ is any curve that starts at the boundary point $z_0$,
moves into $\O$, stays away from $w$, and terminates at $z$.
(Since $z$ is close to the boundary, we don't have to worry about
avoiding the points in $Z$.)

Now the change of variables formula shows that
$$G(z,w)=\text{Re }\alpha(f(z)),$$
where, writing $u=f(z)$,
$$\alpha(u)= \int_{f(\Gamma)} 2A(\zeta)\ d\zeta.$$
Note that $\alpha$ is a local antiderivative of the algebraic function
$2A$ obtained by integrating along a curve that starts at the point
$1$ in the unit circle and moves into the unit disc and terminates
at $u$.  Now, the formula $G(z,w)=\text{Re }\alpha(f(z))$
can be extended to all of $\O$ by means of analytic continuation.
To be more precise, using the notation set up above, if $z$ is in
$\mathcal O_k$, we may define $\alpha$ at $u=f(z)$ as follows.
Let $\gamma$ denote any curve that starts at $1$ and enters the unit
disc and travels to $u$ in $D_1(0)\setminus L$ avoiding $f(w)$.  Define
$$\alpha(u)=\int_\gamma 2A_k(\zeta)\ d\zeta.$$
We may now state that $\text{Re\,}\alpha(f(z))$ is well defined on
$\mathcal O_k$ and that these functions extend continuously up to
$f^{-1}(L)$ minus $w$ and match up there to represent the harmonic function
$G(z,w)$ of $z$ on $\O\setminus \{w\}$.  (Since the functions match
up along $f^{-1}(L)\setminus Z$ minus $w$, we may define a multivalued
version of $\alpha(u)$ globally by integrating a single germ of $2A$ at
$1$ (correctly chosen) along a curve that starts at $1$ and avoids $w$
and the finitely many points in~$W$, where it is understood that we are
integrating the analytic continuation of $2A$ along the curve as we go.)

We remark here that the harmonic measure functions $\omega_j$ can
be dealt with in exactly the same way as we have treated the Green's
function in this section.  Indeed, the same reasoning
that yielded equation~(\ref{eqnC}) can be applied to $\omega_j$ to
obtain
$$F_j'(z)T(z)=-\overline{F_j'(z)}\, \overline{T(z)}$$
on $b\O$, and we may proceed as above to see that $F_j'(z)=f'(z)A(f(z))$
where $A$ is algebraic.  In this way, we prove that $\omega_j$ is
equal to $\text{Re\,}\beta_j(f(z))$ where $\beta_j$ is a Green antiderivative
of type III.

The hard part in what follows is to see how the Green antiderivative
$\alpha$ in the formula for $G(z,w)$ obtained in this section varies
as we allow $w$ to vary.

\section{Proof of Theorem~\ref{thm2}}
\label{sec3}

Assume that $\O$ is a bounded domain in the plane bounded
by $n$ non-intersecting real analytic curves.  We shall need to use some
formulas proved in \cite{B2} that
relate the Poisson kernel to the Szeg\H o kernel $S(z,w)$ and the Garabedian
kernel $L(z,w)$.  Before we write the formulas, we recall some basic facts
about the Szeg\H o and Garabedian kernels on a
domain with real analytic boundary (proofs of which can be found in \cite{B1}).
The kernel $S(z,w)$ extends holomorphically past the boundary in $z$ for each
fixed $w$ in $\O$.  It extends meromorphically past the boundary in $z$ for
each fixed $w$ in $b\O$; in fact, it extends holomorphically past $b\O\setminus \{w\}$
and has only a simple pole at the point $w$.  Furthermore $S(z,w)\ne0$ if
$z\in b\O$ and $w\in\O$.  If $w\in b\O$, then $S(z,w)$ has exactly $n-1$
simple zeroes in $z$, one on each boundary curve different from
the one containing the point $w$.  The kernel $L(z,w)$ has a simple pole in
$z$ at the point $w\in\O$.  It extends holomorphically past the boundary in $z$
for each fixed $w$ in $\O$.  It extends meromorphically past the boundary in $z$
for each fixed $w$ in $b\O$; in fact, it extends holomorphically past $b\O\setminus \{w\}$
and has only a simple pole at the point $w$.  Furthermore $L(z,w)\ne0$ if
$z,w\in \O$ with $z\ne w$.  If $w\in b\O$, then $L(z,w)$ has exactly $n-1$
simple zeroes in $z$, one on each boundary curve different from
the one containing the point $w$ (and these zeroes agree with those of the
Szeg\H o kernel).  Finally, $S(z,w)$ is in
$C^\infty$ of $\Obar\times\Obar$ minus the boundary diagonal
$\{(z,z): z\in b\O\}$ and $L(z,w)$ is in
$C^\infty$ of $\Obar\times\Obar$ minus the diagonal
$\{(z,z): z\in \Obar\}$.

It is proved in \cite{B2} that
\begin{equation}
\label{eqnE}
G_z(z,w)=\pi\frac{S(z,w)L(z,w)}{S(w,w)}+i\pi\sum_{j=1}^{n-1}(\omega_j(w)-\lambda_j(w))u_j(z),
\end{equation}
where the functions $\omega_j$, $\lambda_j$ and $u_j$ were defined in \S\ref{sec1}.
For an alternate and more natural proof of this identity, see \cite[p.~10-12]{B2a}
(where it must be noted that the indices $j$ and $k$ are one and the same).
It is also proved in \cite{B2} that it possible to choose a point $a$ in $\O$
so that the $n-1$ zeroes of $S(z,a)$ in the $z$ variable are distinct and simple.
Choose such a point $a$ and let $a_1,a_2,\dots,a_{n-1}$ denote these simple
zeroes.  For convenience, let $a_0$ denote $a$.  It is proved in \cite{B2}
that the Szeg\H o kernel and Garabedian kernels can be expressed via
\begin{equation}
\label{eqnF}
S(z,w)=\frac{1}{1-\overline{f(w)}\,f(z)}\left(c_{00}S(z,a)\overline{S(w,a)}+
\sum_{j,k=1}^{n-1}c_{jk}S(z,a_j)\overline{S(w,a_k)}\right)
\end{equation}
and
\begin{equation}
\label{eqnG}
L(z,w)=\frac{f(w)}{f(z)-f(w)}\left(c_{00}S(z,a)L(w,a)+
\sum_{j,k=1}^{n-1}\bar c_{jk}S(z,a_j)L(w,a_k)\right).
\end{equation}
To shorten some expressions in what follows, we define coefficients
$c_{0j}=0$ and $c_{j0}=0$ when $j\ne0$.  Let $f$ denote the Ahlfors map
$f_a$ associated to the point $a$.

We now take the principal term
$$X(z,w) := \pi\frac{S(z,w)L(z,w)}{S(w,w)}$$
from equation~(\ref{eqnE}) and replace the Szeg\H o and Garabedian
kernels by the expressions in equations~(\ref{eqnF}) and~(\ref{eqnG}).
Expand this large expression to see that $X(z,w)$ is a linear combination
of rational functions of $f(z)$, $f(w)$, and $\overline{f(w)}$ times
terms of the form
\begin{equation}
\label{eqnH}
\frac{ S(z,a_p)S(z,a_q)\overline{S(w,a_r)}L(w,a_s) }
{ \sum_{j,k=0}^{n-1}c_{jk}S(w,a_j)\overline{S(w,a_k)} }.
\end{equation}
We shall next use the fact that certain combinations of functions extend
to the double as meromorphic functions (in a manner similar to what
was done in \cite{B3a} to study the Carath\'eodory metric and the Poisson
kernel).  In particular
\begin{equation}
\label{eqnI}
\frac{ S(z,a_p)S(z,a_q) }{f'(z)}
\end{equation}
and functions of the form $S(w,a_r)/S(w,a)$ and $L(w,a_s)/S(w,a)$
all extend meromorphically to the double.  To see that
$S(w,a_r)/S(w,a)$ extends to the double, use the standard identity
$$\overline{S(z,w)}=\frac{1}{i}L(z,w)T(z),$$
which holds for $z\in b\O$ and $w\in\O$, to see that $S(w,a_r)/S(w,a)$
is equal to the conjugate of $L(w,a_r)/L(w,a)$ on the boundary.  Hence,
it extends to the double as a meromorphic function.  Similarly,
$L(w,a_s)/S(w,a)$ is equal to the conjugate of $S(w,a_s)/L(w,a)$
on the boundary, and so it extends to the double.  To see that
$S(z,a_p)S(z,a_q)/f'(z)$ extends, we use identity~(\ref{eqnD}) and
the fact that $\overline{f(z)}=1/f(z)$ on $b\O$.  Indeed, (\ref{eqnI}) is
equal to the conjugate of $f(z)^2 L(z,a_p)L(z,a_q)/f'(z)$ on the boundary,
and so it extends meromorphically to the double. 
We may now multiply (\ref{eqnH}) by unity three times, once in the form
$f'(z)/f'(z)$ and in the form $S(w,a)/S(w,a)$ and the conjugate of
$S(w,a)/S(w,a)$.  Next, we distribute the terms $1/f'(z)$, $1/S(w,a)$ and 
the conjugate of $1/S(w,a)$ into the expression in strategic places in
order to rewrite (\ref{eqnH}) in the form
\begin{equation}
\label{eqnJ}
f'(z)\ \frac{
\frac{S(z,a_p)S(z,a_q)}{f'(z)}\ \frac{\overline{S(w,a_r)}}{\overline{S(w,a)}}
\ \frac{L(w,a_s)}{S(w,a)}
}
{
\sum_{j,k=0}^{n-1}c_{jk}\frac{S(w,a_j)}{S(w,a)}
\ \frac{\overline{S(w,a_k)}}{\overline{S(w,a)}}
}.
\end{equation}
Every quotient in this last expression extends to the double as a
meromorphic or antimeromorphic function.  As mentioned earlier, we proved
in \cite{B3} that the field of meromorphic functions on the double is
generated by $f$ and another Ahlfors map $F$.  Hence, we have shown that
$X(z,w)$ is equal to $f'(z)$ times a rational function of $f(z)$, $F(z)$,
$f(w)$, $F(w)$, $\overline{f(w)}$, and $\overline{F(w)}$.  Since $f$
and $F$ are algebraically dependent, $F$ is an algebraic function of $f$
and we conclude that
$$X(z,w)=f'(z)A(f(z),f(w),\overline{f(w)}\,),$$
where $A(\zeta,u,v)$ is an algebraic function of three variables.
Let $z_0$ denote the point in the outer boundary of $\O$ such that $f(z_0)=1$.
We may now repeat the argument in \S\ref{sec2} and integrate equation~(\ref{eqnE})
in the first variable along a curve $\Gamma$ that starts at  $z_0$, avoids $w$
and the finitely many points in the branch locus of $f$, and terminates at $z$.
Note that, since $2(\dee/\dee z)\mu_j=iu_j$, and since $\mu_j$ is real valued and
vanishes on the outer boundary, it follows that
$\text{Re\,}\int_\Gamma i u_j dz= \mu_j(z)$, and we obtain
the formula for the Green's function in the statement of Theorem~\ref{thm2}.
This completes the proof in case $\O$ has real analytic boundary.  In the more
general case of the statement of the theorem, we note that all the elements
in the formula transform under a biholomorphic mapping to a domain with real
analytic boundary and so the formula holds in this more general setting.

\section{Proof of Theorem~\ref{thm1}}
\label{sec4}
As explained in \S\ref{sec2}, we may assume that $\O$ is a bounded domain bounded
by $n$ non-intersecting real analytic Jordan curves.  We proved in \S\ref{sec3} that
\begin{equation}
\label{eqnK}
G_z(z,w)=f'(z)A(f(z),f(w),\overline{f(w)}\,)+
i\pi\sum_{j=1}^{n-1}(\omega_j(w)-\lambda_j(w))u_j(z),
\end{equation}
where $A$ is algebraic.  Recall that
$$u_j=\sum_{k=1}^{n-1}\sigma_{jk} F_k',$$
where the coefficients are such that $i\sigma_{jk}$ is real.
Hence, the sum
$$i\pi\sum_{j=1}^{n-1}(\omega_j(w)-\lambda_j(w))u_j(z)$$
can be rearranged to appear in the form
$$\sum_{j=1}^{n-1}v_j(w)F_j'(z),$$
where $v_j$ is a real valued function that is a linear combination over
the real numbers of the real valued functions $\omega_j-\lambda_j$.
Let $z_0$ be the point on the outer boundary such that $f(z_0)=1$,
and let $\nu_k$ denote a curve that starts at $z_0$, enters
$\O$ and avoids $f(w)$ and the branch points of $f$ and terminates at the
boundary point $z_k$ on $\gamma_k$ such that $f(z_k)=1$.
Integrate equation~(\ref{eqnK}) with respect to $dz$ along $\nu_k$
and take the real part, noting that $2\text{Re\,}\int_{\nu_k}F_j'(z)\,dz =
\omega_j(z_k)-\omega_j(z_0)=1$
if $j=k$ and zero otherwise, to obtain
$$0=2\text{Re }\int_{\nu_k} f'(z)A(f(z),f(w),\overline{f(w)}\,)\ dz
+v_k(w).$$
Finally, we may use the change of variables formula as in \S\ref{sec2} to
see that
$$
v_k(w)
=-\text{Re }\int_{f(\nu_k)} 2A(\zeta,f(w),\overline{f(w)}\,)\ d\zeta,$$
where $f(\nu_k)$ is a curve that starts at $1$, enters the unit disc,
avoids the point $f(w)$ and the finitely many points in the image of
the branch locus of $f$, and terminates back at $1$.  This shows that
$v_k$ is equal to $2\alpha_k(f(w),\overline{f(w)}\,)$ where $\alpha_k$
is a Green antiderivative of type~II.

Now, when we integrate
$$G_z(z,w)=f'(z)A(f(z),f(w),\overline{f(w)}\,)+
\sum_{j=1}^{n-1}\alpha_j(f(w),\overline{f(w)}\,)2F_j'(z)$$
with respect to $z$ along a curve $\Gamma$ as we did in \S\ref{sec3}, we
obtain
$$G(z,w)=\text{Re\,}\alpha(f(z),f(w),\overline{f(w)}\,)+
\sum_{j=1}^{n-1}\alpha_j(f(w),\overline{f(w)}\,)\omega_j(z).$$
We showed in \S2 that $\omega_j(z)$ is equal to $\beta_j(f(z))$
where $\beta_j$ is a Green antiderivative of type~III.  This
completes the proof of Theorem~\ref{thm1}.

\section{Proof of Theorem~\ref{thm3}}
\label{sec5}
We continue to assume that $\O$ is a bounded domain bounded by $n$
non-intersecting real analytic Jordan curves.
Integrate equation~(\ref{eqnK})
with respect to $dz$ around one of the inner boundary curves,
say $\gamma_k$.  Note that, because the Poisson kernel is
given by $P(w,z)=(-i/\pi)G_z(z,w)T(z)$, it follows that
$$\int_{\gamma_k}G_z(z,w)\ dz=i\pi\omega_k(w).$$
Also, recall that
$\int_{\gamma_k}u_k(z)\ dz=1$ and
$\int_{\gamma_k}u_j(z)\ dz=0$ if $j\ne k$.
Hence,
$$ \omega_k(w)=
\frac{1}{i\pi}
\int_{\gamma_k} f'(z)A(f(z),f(w),\overline{f(w)}\,)\ dz +
(\omega_k(w)-\lambda_k(w)),$$
and consequently
$$
\lambda_k(w)=
\frac{1}{i\pi}
\int_{\gamma_k} f'(z)A(f(z),f(w),\overline{f(w)}\,)\ dz.
$$
We may now use the change of variables formula to see that
$$
\lambda_k(w)=
\frac{1}{i\pi}
\int_{C_1} A(\zeta,f(w),\overline{f(w)}\,)\ d\zeta,
$$
where $C_1$ denotes the unit circle taken in the standard sense.
This shows that
$\lambda_k(w)$ is equal to $\text{Re }\psi_k(f(w),\overline{f(w)}\,)$
where $\psi_k$ is a Green antiderivative of type~II.  (In fact, the
real part can be dropped here because $\psi_k(f(w),\overline{f(w)}\,)$
is real valued on $\O$.)

We proved in $\S2$ that $\omega_j(z)$ is equal to $\beta_j(f(z))$
where $\beta_j$ is a Green antiderivative of type~III.  Since $\mu_k$
is a linear combination of the functions $\omega_j$, $j=1,\dots,n-1$,
it is equal to $\tau_k(f(z))$ where $\tau_k$ is a Green antiderivative
of type~III.  This completes the proof of Theorem~\ref{thm3}.

\section{Extension of the Green's function}
\label{sec6}

A fascinating consequence of the formulas in this paper is that
the Green's function associated to a domain with an algebraic
Ahlfors map can be harmonically continued in both variables along
curves in the complex plane that avoid finitely many points.  At those
finitely many points, the singular behavior is limited by the
nature of an antiderivative of an algebraic function composed
with an algebraic function.  Such extension behavior was intimated by
results of Ebenfelt \cite{E} and Khavinson and Shapiro \cite{KhS} on
extension of solutions to the Dirichlet problem.  We remark that Aharonov
and Shapiro \cite{AS} proved that the Ahlfors map associated to
a quadrature domain in the plane must be algebraic, and Gustafsson
\cite{G} proved that quadrature domains are dense among all
bounded finitely connected domains in the plane bounded by Jordan
curves.  In fact, it was proved in \cite{B4} (see also \cite{B5}) that,
given a bounded finitely connected domain bounded by $C^\infty$
smooth Jordan curves, it is possible to find a quadrature domain
that is arbitrarily $C^\infty$ close by which is biholomorphic
via a mapping which is arbitrarily $C^\infty$ close to the identity
map.  Thus, it is possible to make very subtle changes in a smooth
domain so that the Green's function continues harmonically to the whole
complex plane minus finitely many points.

In general, we have shown that the Green's function extends harmonically
to the set where the Ahlfors map extends holomorphically (minus finitely
many points).  We would be very curious to know what the Green's
function is for the domain $\{z:|z+(1/z)|<r\}$, where $r$ is a real
constant greater than~$2$.  Could it be an elementary function in the
sense of Liouville?  It may be the simplest Green's function
associated to a multiply connected domain, but we are willing
to bet that no Green's function associated to a multiply connected
domain can be an elementary function.

\end{document}